\documentclass[a4paper, 10pt, conference]{ieeeconf}      

\IEEEoverridecommandlockouts                              

\usepackage{graphicx}      
\usepackage{amsmath} 
\usepackage{amssymb}  
\usepackage{epsfig}
\usepackage{graphicx}
\usepackage[hang]{caption2}
\usepackage[normalsize,it,hang]{subfigure}
\usepackage[latin1]{inputenc}
\usepackage[greek,english]{babel}
\usepackage{longtable}
\usepackage{booktabs}

\newtheorem{remark}{Remark}[section]

\title{\LARGE \bf
Stability margins and model-free control: A first look
}

\author{Michel FLIESS$^{a,c}$ and C\'edric JOIN$^{b,c,d}$ 
\thanks{$^{a}$ LIX (CNRS, UMR 7161), \'Ecole polytechnique, 91128 Palaiseau, France.
        {\tt\small Michel.Fliess@polytechnique.edu}}%
\thanks{$^{b}$ CRAN (CNRS, UMR 7039), Universit\'e de Lorraine, BP 239, 54506 Vand{\oe}uvre-l\`es-Nancy, France.
 \newline {\tt\small cedric.join@univ-lorraine.fr}}%
 \thanks{$^{c}$ AL.I.E.N. (ALgèbre pour Identification \& Estimation Numériques), 24-30 rue Lionnois, BP 60120, 54003 Nancy, France.
 \newline       {\tt\small \{michel.fliess, cedric.join\}@alien-sas.com}}%
 \thanks{$^{d}$ Projet NON-A, INRIA Lille -- Nord-Europe, France.} 
 }

\begin{document}

\maketitle
\thispagestyle{empty}
\pagestyle{empty}

\begin{abstract}
We show that the open-loop transfer functions and the stability margins may be defined within the recent model-free control setting.  Several convincing computer experiments are 
presented including one which studies the robustness with respect to delays. 
\end{abstract}

\section{Introduction}

 \emph{Stability margins} are basic ingredients of control theory. They are widely taught (see, {\it e.g.}, \cite{murray,bourles,follinger,Franklin,kuo,lausanne}, and the references therein) and are quite 
 often utilized in industry in order to check the control design of plants, or, more exactly, of their mathematical models. The importance of this topic is highlighted by the following fact:
 the literature on theoretical advances and on the connections with many case-studies 
 contains several thousands of publications! This communication relates stability margins to the recent \emph{model-free} control and the corresponding \emph{intelligent} PIDs
 \cite{ijc13}, which were already illustrated by many concrete and varied applications (see, \textit{e.g.}, the numerous references in \cite{ijc13}, and, during the last months, \cite{compiegne,lafont,pol,menhour,michel1,michel2,tarbes,wang,lyon}).
 \begin{remark}
Let us emphasize that our model-free control design and the corresponding intelligent controllers are most easily implementable (see \cite{ijc13,nice}).
\end{remark}
Our aims are the following ones:
 \begin{enumerate}
\item Practitioners of stability margins and other frequency techniques will recognize that their expertise still makes sense within model-free control.
\item The influence of delays in model-free control is analyzed for the first time.
\end{enumerate}

Let us briefly explain our viewpoint. Take a monovariable system which is governed  by unknown equations. Consider the \emph{ultra-local} model \cite{ijc13}
\begin{equation}
\boxed{\dot{y} = F + \alpha u } \label{ul1}
\end{equation} 
where 
\begin{itemize}
\item $u$ and $y$ are respectively the input and output variables,
\item $F$ subsumes the unknown parts, including the perturbations,
\item $\alpha$ is a constant parameter which is chosen by the engineer in such a way that $\alpha u$ and $\dot{y}$ are of the same magnitude.
\end{itemize}
Remember that Equation \eqref{ul1} applies not only to systems with lumped parameters, \textit{i.e.}, to systems which are described by ordinary differential equations 
of any order, but also to systems with distributed parameters, \textit{i.e.}, to partial differential equations (see, \textit{e.g.}, \cite{edf}). 
Close the loop with
 \begin{equation}\label{feedbackintro}
 u = \frac{- F_{\text{est}} + \dot{y}^\star + \mathfrak{F}}{\alpha} 
 \end{equation}
such that
\begin{itemize}
\item   $F_{\text{est}}$ is a realtime estimate of $F $ (see Section \ref{F}),
\item ${y}^\star $ is a reference trajectory,
\item the closed loop system is
\begin{equation}\label{fd1}
\dot{e} +  \mathfrak{F} = F_{\text{est}} - F 
\end{equation} 
where $e =y^\star- y$ is the tracking error,
\item $\mathfrak{F}$ is either a proportional controller
\begin{equation}\label{fdp}
\mathfrak{F} = K_P e
\end{equation} 
or, sometimes, a proportional-integral controller
\begin{equation}\label{fdpi}
\mathfrak{F} = K_P e + K_I \int e 
\end{equation} 
such that
\begin{equation}\label{fd2}
\dot{e} + \mathfrak{F} = 0
\end{equation} 
exhibits the desired asymptotic stability. For instance $K_P$ in Equation \eqref{fdp} should be positive.
\end{itemize}
Equations \eqref{fdp}-\eqref{fd2} and  \eqref{fdpi}-\eqref{fd2}  yield the usual open-loop transfer functions 
\begin{equation}\label{trans1}
\boxed{T_{1OLP} = \frac{K_P}{s}}
\end{equation} 
and
\begin{equation}\label{trans2}
\boxed{T_{1OLPI} =  \frac{1}{s} \times (K_P + \frac{K_I}{s}) = \frac{K_P}{s} +\frac{K_I}{s^2}}
\end{equation}
Their \emph{gain} and \emph{phase margins}  are by definitions those of the systems defined by Equations \eqref{trans1} and \eqref{trans2}. 
Note that $F_{\text{est}} - F$ in Equation \eqref{fd1} should be viewed as an additive disturbance.

Our paper is organized as follows. Basics of model-free control are briefly revisited in Section \ref{mfc}. Section \ref{SMFC} computes some open loop functions for iPID's, iPD's, iPIs, iPs as well as the corresponding 
stability margins. Several computer experiments are examined in Section \ref{ILL}, including the robustness with respect to delays. Concluding remarks are developed in Section \ref{conclusion}.

\section{Model-free control: A short review\protect\footnote{See \cite{ijc13} for more details.}}\label{mfc}

\subsection{The ultra-local model}
Introduce the \emph{ultra-local model}
\begin{equation}
\boxed{y^{(\nu)} = F + \alpha u} \label{ultralocal}
\end{equation}
where
\begin{itemize}
\item the order $\nu \geq 0$ of derivation is a non-negative
integer which is selected by the practitioner,\footnote{The existing
examples show that $\nu$ may always be chosen quite low,
\textit{i.e.}, $1$, or $2$. Most of the times $\nu = 1$. The only concrete example until now with $\nu = 2$ is provided by the 
magnetic bearing \cite{compiegne}, where the friction is negligible (see the explanation in \cite{ijc13}). }
\item $\alpha \in \mathbb{R}$ is chosen by the practitioner such that $\alpha u$ and
$y^{(\nu)}$ are of the same magnitude,
\item $F$ represents the unknown structure of the control system as well as the perturbations. 
\end{itemize}

\subsection{Intelligent controllers}
\subsubsection{Generalities}
Close the loop with respect to Equation \eqref{ultralocal} via the
\emph{intelligent controller}
\begin{equation}\label{ic}
u = \frac{- F_{\text{est}} + y^{\ast (\nu)} + \mathfrak{F}(e)}{\alpha}
\end{equation}
where
\begin{itemize}
\item $F_{\text{est}} $ is a realtime estimate of $F$,
\item $y^\ast$ is the output reference trajectory,
\item $e = y^\ast - y$ is the tracking error,
\item $\mathfrak{F}(e)$ is a functional of $e$ such that the closed-loop system
\begin{equation}\label{fc}
e^{(\nu)} + \mathfrak{F}(e) =F_{\text{est}} - F 
\end{equation}
exhibits a desired behavior. If, in particular, the estimation is perfect, \textit{i.e.}, $F_{\text{est}}  = F$, then 
\begin{equation}\label{fco}
e^{(\nu)} + \mathfrak{F}(e) = 0
\end{equation}
should be asymptotically stable, \textit{i.e.}, $\lim_{t\rightarrow +\infty} e(t) = 0$.
\end{itemize}
\subsubsection{Intelligent PIDs}
If $\nu = 2$ in Equation \eqref{ultralocal}, \textit{i.e.},
\begin{equation}
\boxed{\ddot{y} = F + \alpha u} \label{ord2}
\end{equation}
Close the loop via the \emph{intelligent
proportional-integral-derivative controller}, or \emph{iPID},
\begin{equation}\label{ipid}
\boxed{u =  \frac{- F_{\text{est}} + \ddot{y}^\ast + K_P e + K_I \int e + K_D
\dot{e}}{\alpha}}
\end{equation}
where $K_P$, $K_I$, $K_D$ are the usual tuning gains. Combining
Equations \eqref{ord2} and \eqref{ipid} yields
$$
\ddot{e} + K_D \dot{e} + K_P e + K_I \int e = F_{\text{est}}-F 
$$
$K_I = 0$ in Equation \eqref{ipid} yields the \emph{intelligent
proportional-derivative controller}, or \emph{iPD},
\begin{equation}\label{ipd}
\boxed{u =  \frac{- F_{\text{est}} + \ddot{y}^\ast + K_P e + K_D
\dot{e}}{\alpha}}
\end{equation}
Such an iPD was employed in \cite{compiegne}.

If $\nu = 1$ in Equation \eqref{ultralocal}, we recover Equation \eqref{ul1}.
The loop is closed by the \emph{intelligent proportional-integral
controller}, or \emph{iPI},
\begin{equation*}\label{ipi}
\boxed{u = \frac{- F_{\text{est}}  + \dot{y}^\ast + K_P e  + K_I \int e}{\alpha}}
\end{equation*}
$K_I$ may often be set to $0$. It yields the \emph{intelligent
proportional controller}, or \emph{iP},
\begin{equation*}\label{ip}
\boxed{u = - \frac{- F_{\text{est}} + \dot{y}^\ast +  K_P e}{\alpha}}
\end{equation*}

\subsection{Estimation of $F$}\label{F}
Assume that $F$ in Equation \eqref{ultralocal} may be ``well'' approximated by a piecewise constant function $F_{\text{est}} $. According to the algebraic parameter identification 
developed in \cite{sira1,sira2}, rewrite, if $\nu = 1$, Equation \eqref{ul1}  in the operational domain (see, \textit{e.g.}, \cite{yosida}) 
$$s
Y = \frac{\Phi}{s}+\alpha U +y(0)
$$
where $\Phi$ is a constant. We get rid of the initial condition $y(0)$ by multiplying both sides on the left by $\frac{d}{ds}$:
$$
Y + s\frac{dY}{ds}=-\frac{\Phi}{s^2}+\alpha \frac{dU}{ds}
$$
Noise attenuation is achieved by multiplying both sides on the left by $s^{-2}$. It yields in the time domain the realtime estimate
$$
F_{\text{est}}(t)  =-\frac{6}{\tau^3}\int_{t-\tau}^t \left( (\tau -2\delta)y(\delta)+\alpha\delta(\tau -\delta)u(\delta) \right)d\delta
$$
where $\tau > 0$ might be ``small''. 
\begin{remark}
As in our first publications on model-free control, $F_{\text{est}}(t)$ might also be obtained by estimating the noisy derivative of $y$ (see \cite{nl}, and \cite{mboup}, \cite{liu}). 
\end{remark}

\section{Open-loop transfer functions}\label{SMFC}
\subsection{Definitions}
Assume that in Equation \eqref{ic} $\mathfrak{F}$ may be defined by a transfer function $T_\mathfrak{F}$. Then Equation \eqref{fco} yields the 
transfer function 
\begin{equation}\label{ol}
T_{\nu OL} = \frac{T_\mathfrak{F}}{s^\nu}
\end{equation}
which is called the \emph{open-loop transfer function} of the system defined by Equations \eqref{ultralocal} and \eqref{ic}. If $\nu = 2$, and 
with an iPID, the open-loop transfer function \eqref{ol} of the system defined by Equations \eqref{ord2}  and \eqref{ipid} becomes 
\begin{equation}\label{2ol}
\boxed{T_{2OLPID} = \frac{K_P}{s^2} + \frac{K_I}{s^3} + \frac{K_D}{s}  }
\end{equation}
It becomes for an iPD:
\begin{equation}\label{2olpd}
\boxed{T_{2OLPD} = \frac{K_P}{s^2} + \frac{K_D}{s}  }
\end{equation}
If $\nu = 1$, and 
with an iPI, the open-loop transfer function of the system defined by Equations \eqref{ul1}  and \eqref{fdp}  or \eqref{fdpi}, the corresponding open-loop transfer functions
were already given by Equations \eqref{trans1} and \eqref{trans2}.

\begin{remark}\label{21}
Notice that $T_{1OLPI} $ and $T_{2OLPD} $ are expressed by Formulae \eqref{trans2} and \eqref{2olpd} , which are identical if 
we exchange $K_P$, $K_I$ with $K_I$, $K_D$.
\end{remark}

\subsection{Stability margins}

\subsubsection{iP}\label{IP}
Setting $s = j \omega$ in Equation \eqref {trans1}, where 
\begin{itemize}
\item $\omega$ is a non-negative real number, 
\item $j = \sqrt{- 1}$,
\end{itemize}
gives $T_{1OLP}(jw)=\frac{K_P}{j\omega} = -j \frac{K_P}{\omega} $. Since $K_P > 0$ and $\omega \geq 0$, we obtain the following margins:
 $$\text{PhaseMargin}_{1OLP} = 90°$$
and
$$\text{GainMargin}_{1OLP}=+\infty$$

\subsubsection{iPI}\label{MiP}
Setting as above $s = j \omega$ in Equation \eqref {trans2} yields a complex quantity where the imaginary part is $- j \frac{K_P}{\omega}$. Therefore
$$\text{GainMargin}_{1OLPI} = +\infty$$
and 
 $$\text{PhaseMargin}_{1OLPI}=\tan^{-1}\left(\frac{K_P \omega_{m}}{K_I}\right)$$
 where 
 $$\omega_{m}=\sqrt{\frac{K_P^2+\sqrt{K_P^4+4K_I^2}}{2}}$$ 
 is such that the module of $T_{1OLPI}$ is equal to $1$.
A phase margin of $45°$, for instance, is obtained by setting 
$$\omega_m = \frac{K_I}{K_P}$$
$K_P$ and $K_I$ are then related by the equation
$$
 \frac{K_I}{K_P} = \sqrt{\frac{K_P^2+\sqrt{K_P^4+4K_I^2}}{2}}
 $$

\subsubsection{iPD}
It suffices according to Remark \ref{21}  to replace, in the expressions related to the iPIs in Section \ref{MiP}, 
$K_P$ and $K_I$ respectively by $K_D$ and $K_I$.

\subsubsection{iPID}
It follows from Equation \eqref{2ol} that the stability margins necessitates here the famous Cardano formulae which give the roots of third degree 
algebraic equations (see, \textit{e.g.}, \cite{waerden}). 
A single root is moreover real. Then
$$\text{GainMargin}_{2OLPID} = \frac{K_I}{K_D K_P}$$
and
 $$\text{PhaseMargin}_{2OLPID} = \tan^{-1}\left(\frac{K_D \omega_{m}^2-K_I}{K_Pw_{m}}\right)$$
 where
\begin{equation*}
\omega_{m} = \sqrt{A+\frac{B}{C}+D}
\end{equation*}
$A$, $B$, $C$, $D$ are given by
\begin{equation*}
\begin{split}
A=&\\
&\left(\frac{K_I^2}{2} + \frac{K_P^6}{27} - \frac{K_P^2(2K_IK_P - K_D^2)}{6}\right)^{2/6}\\
&- \left(\frac{K_D^2}{3} + \frac{K_P^4}{9} - \frac{2K_IK_P}{3}\right)^{3/6}\\
&+ \left(\frac{K_I^2}{2} + \frac{K_P^6}{27} - \frac{K_P^2(2K_IK_P - K_D^2)}{6}\right)^{1/3}
\end{split}
\end{equation*}
$$B=\frac{K_D^2}{3} + \frac{K_P^4}{9} - \frac{2K_IK_P}{3}$$
\begin{equation}
\begin{split}
C=&\\
&\left(\frac{K_I^2}{2} + \frac{K_P^6}{27} - \frac{K_P^2(2K_IK_P - K_D^2)}{6}\right)^{2/3}\\
&- \left(\frac{K_D^2}{3} + \frac{K_P^4}{9} - \frac{2K_IK_P}{3}\right)^{1/6}\\
&+ \left(\frac{K_I^2}{2} + \frac{K_P^6}{27} - \frac{K_P^2(2K_IK_P - K_D^2)}{6}\right)^{1/3}
\end{split}
\end{equation}

$$D=\frac{K_P^2}{3}$$

\begin{remark}
Since model-free control encompasses to some extent nonlinear control, the above calculations yield a kind of nonlinear generalization of stability margins (see Section \ref{nla}). Remember that the stability margins for nonlinear
systems have been studied in a number of publications (see, \textit{e.g.}, \cite{sepulchre}).
\end{remark}

\section{Numerical illustrations}\label{ILL}
The equations of the systems considered below are only given for achieving of course computer simulations.

\subsection{A nonlinear academic example}\label{nla}
\subsubsection{Description and control}\label{NLDES}

Consider the stable single-input single-output system
\begin{equation}\label{nl0}
\ddot y+4\dot y+3y = 3\dot uu^2+2u^3
\end{equation}
Our ultra-local model is
\begin{equation}\label{nl1}
\dot y =F+u
\end{equation}
\textit{i.e.}, $\nu=1$, $\alpha=1$ in Equation \eqref{ultralocal}. We employ the iP controller \eqref{fdp} where $K_P=1$.
The gain and phase margins are given by Section \ref{IP}, as in many concrete systems. 
\subsubsection{Some computer experiments}
\begin{figure*}
\begin{center}
\subfigure[Control]{
\resizebox*{5.565cm}{!}{\includegraphics{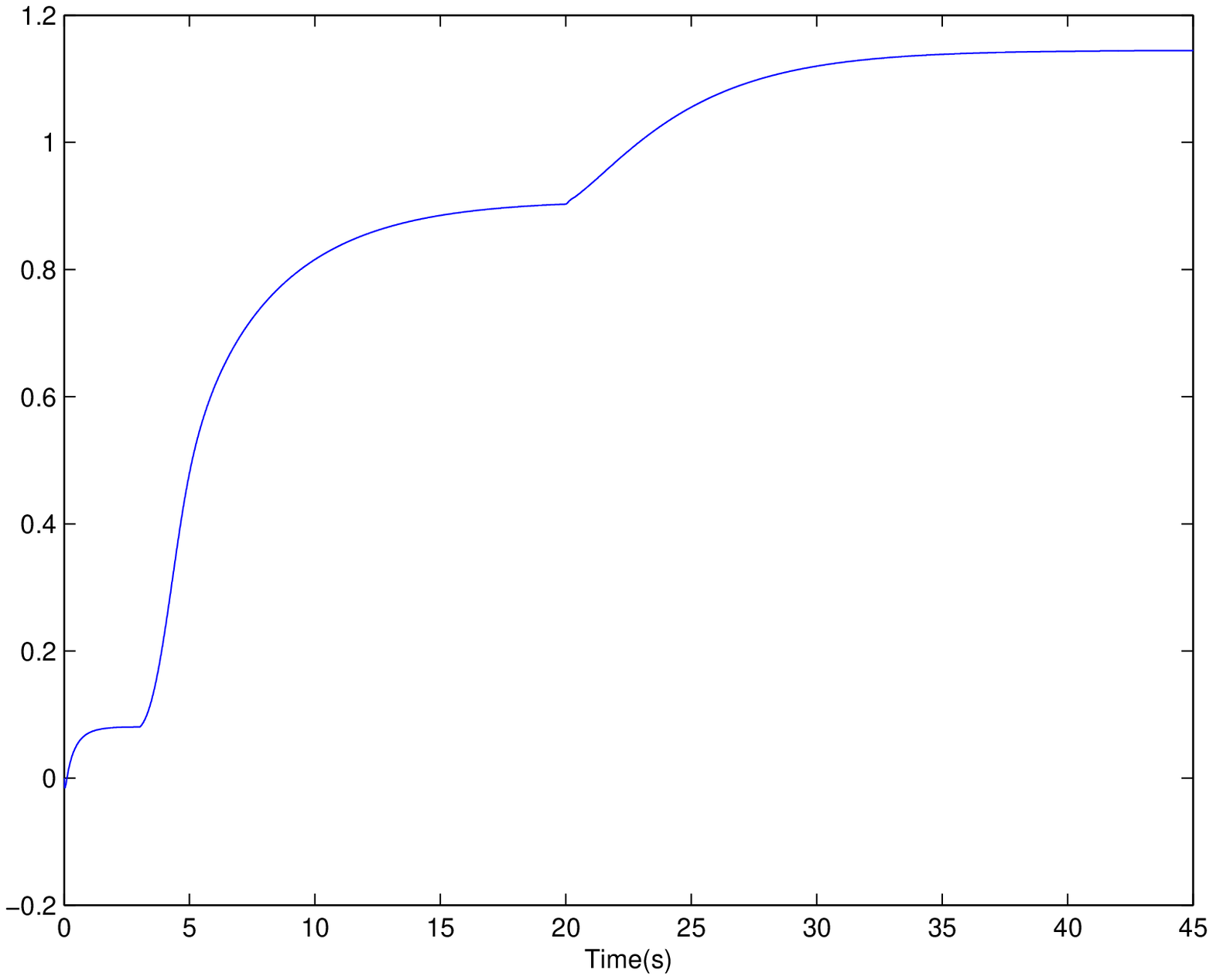}}}%
\subfigure[Setpoint (- .,  black), Reference (- -,  red) and Output (--,  blue)]{
\resizebox*{5.565cm}{!}{\includegraphics{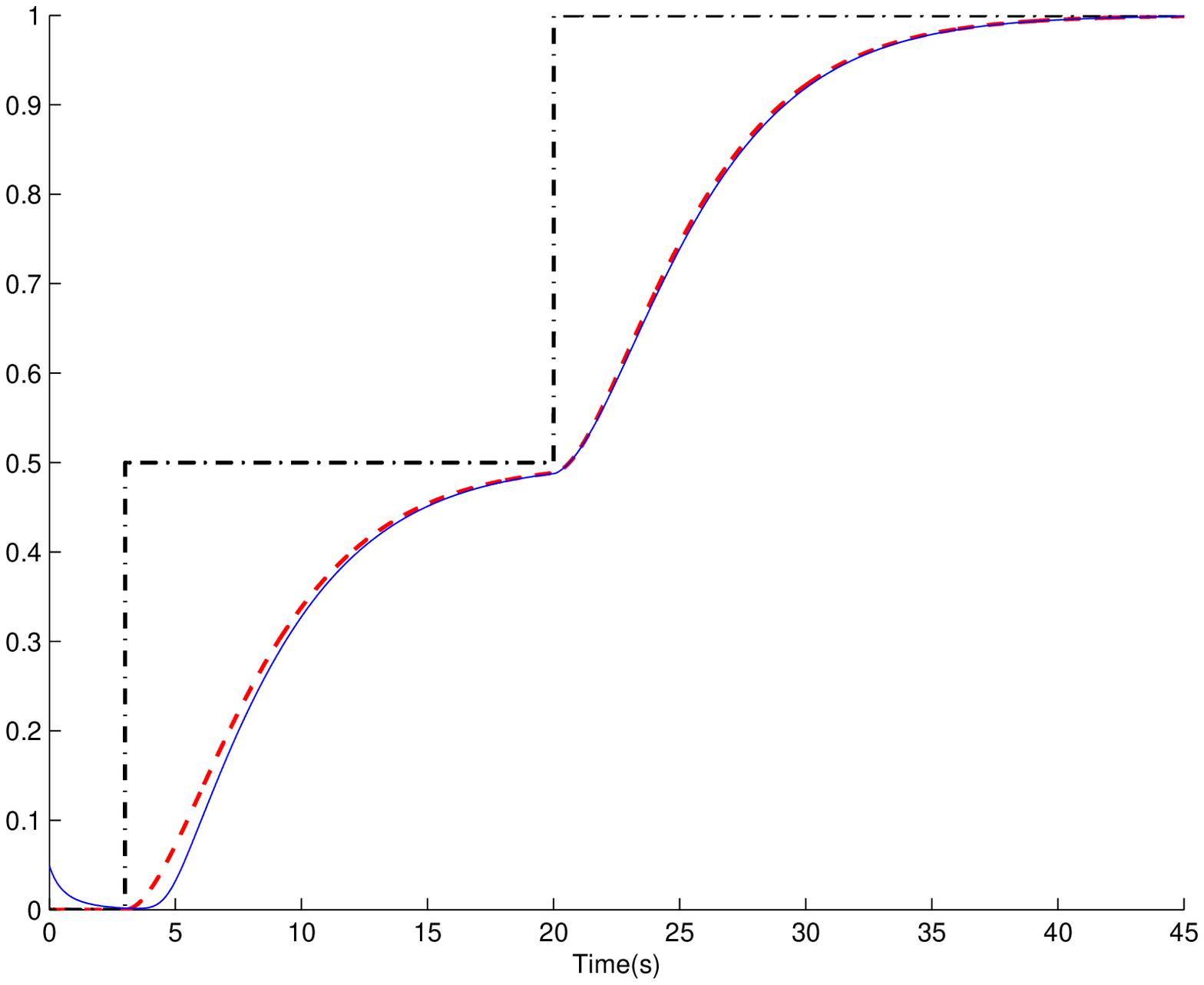}}}
\subfigure[$F$ (- -,  red) and $F_{est}$ (--,  blue)]{
\resizebox*{5.565cm}{!}{\includegraphics{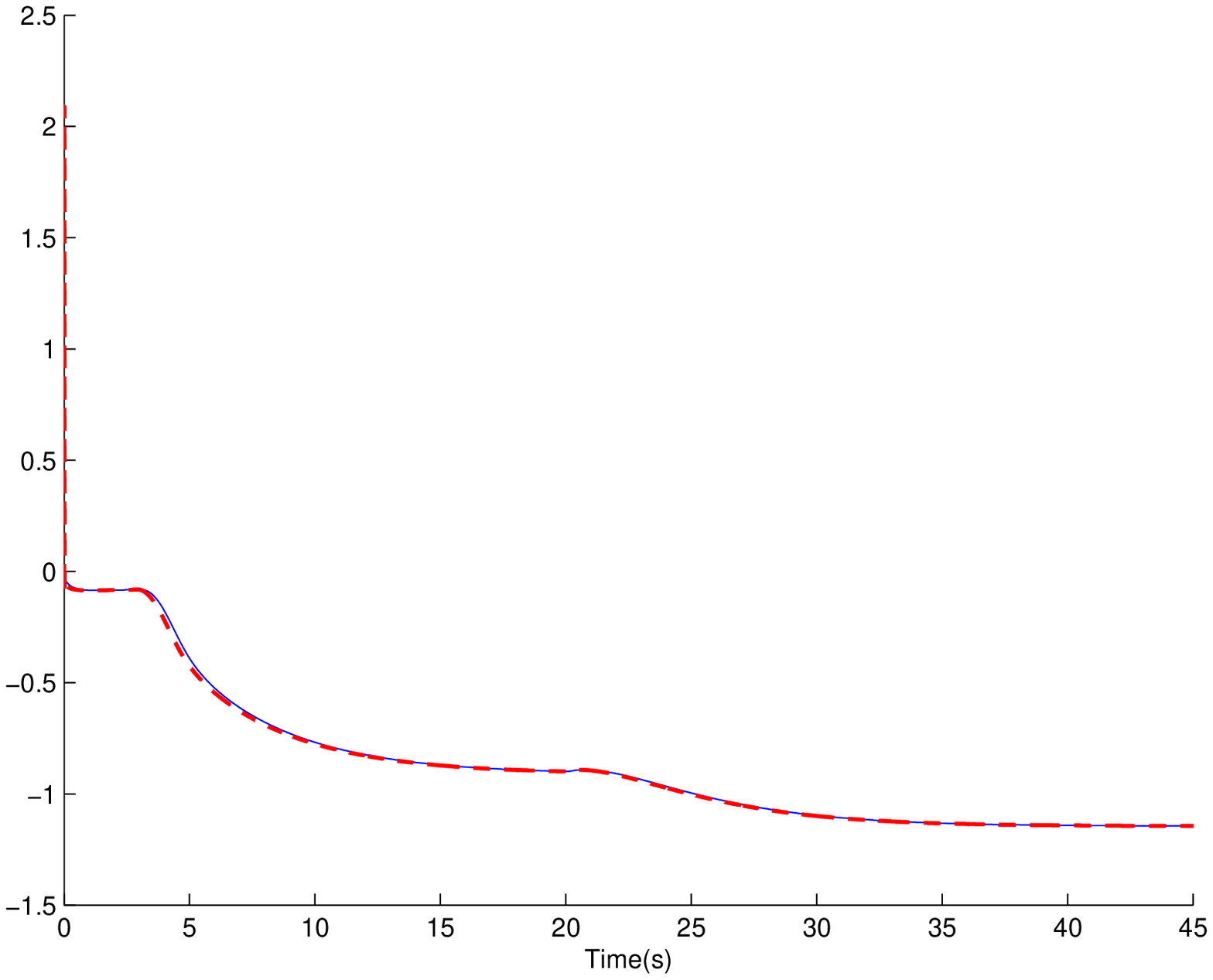}}}%
\caption{$T_{\text{est}}  = 0.01$s}
\label{F11}
\end{center}
\end{figure*}

\begin{figure*}
\begin{center}
\subfigure[Control]{
\resizebox*{5.565cm}{!}{\includegraphics{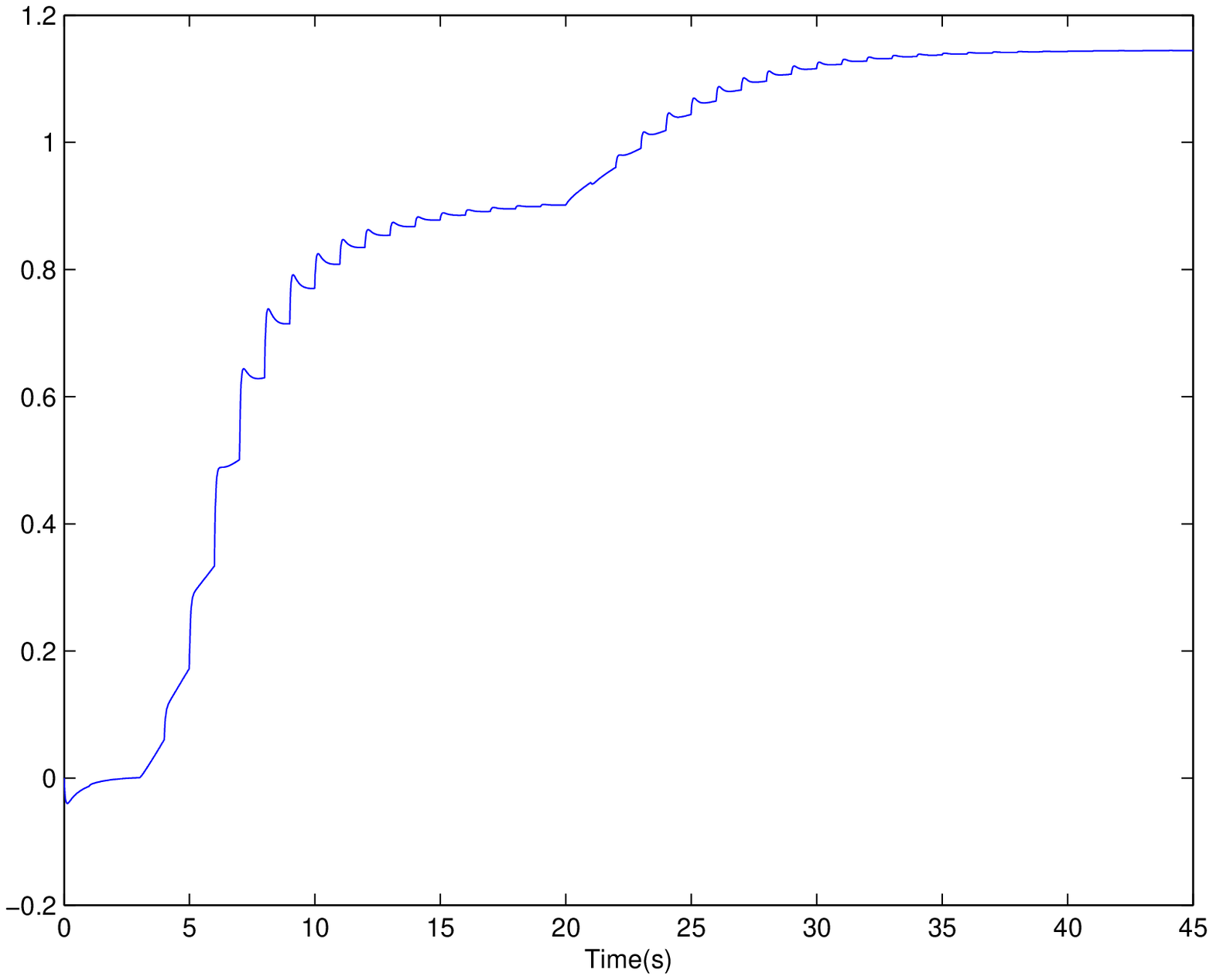}}}%
\subfigure[Setpoint (- .,  black), Reference (- -,  red) and Output (--,  blue)]{
\resizebox*{5.565cm}{!}{\includegraphics{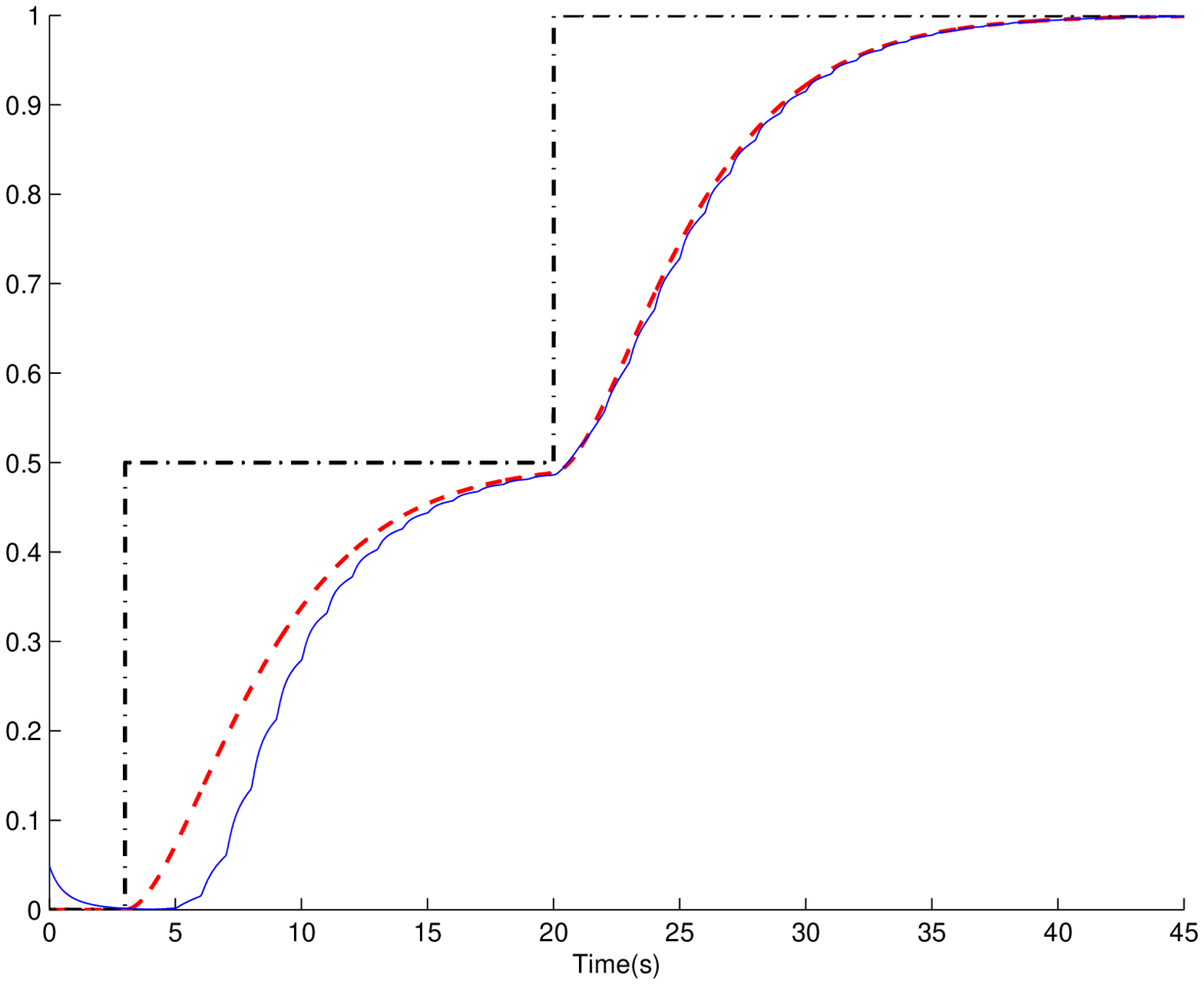}}}
\subfigure[$F$ (- -, red) and $F_{est}$ (--, blue)]{
\resizebox*{5.565cm}{!}{\includegraphics{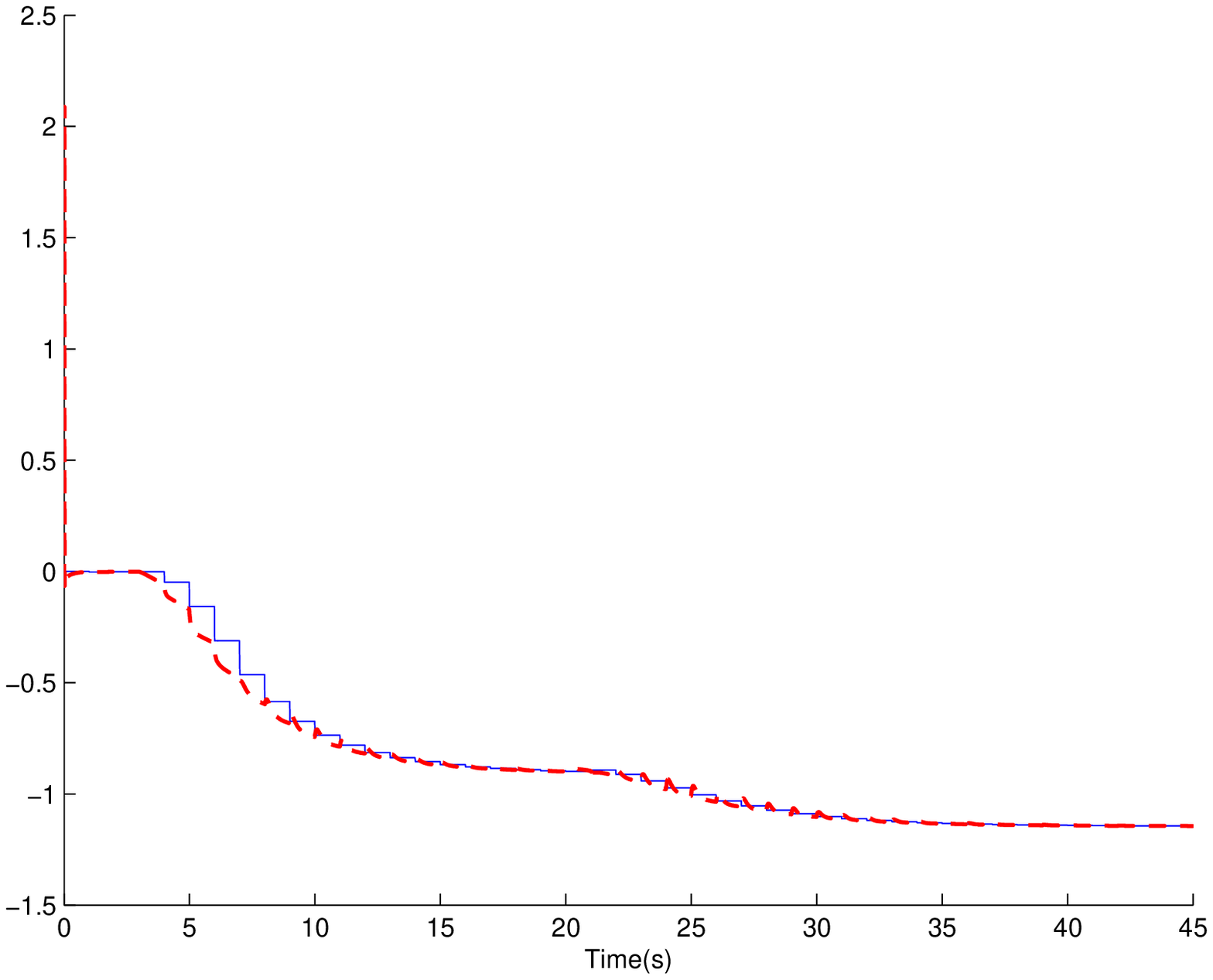}}}%
\caption{$T_{\text{est}}  = 1\text{s}$}
\label{F100}
\end{center}
\end{figure*}

\begin{figure*}
\begin{center}
\subfigure[Control]{
\resizebox*{5.565cm}{!}{\includegraphics{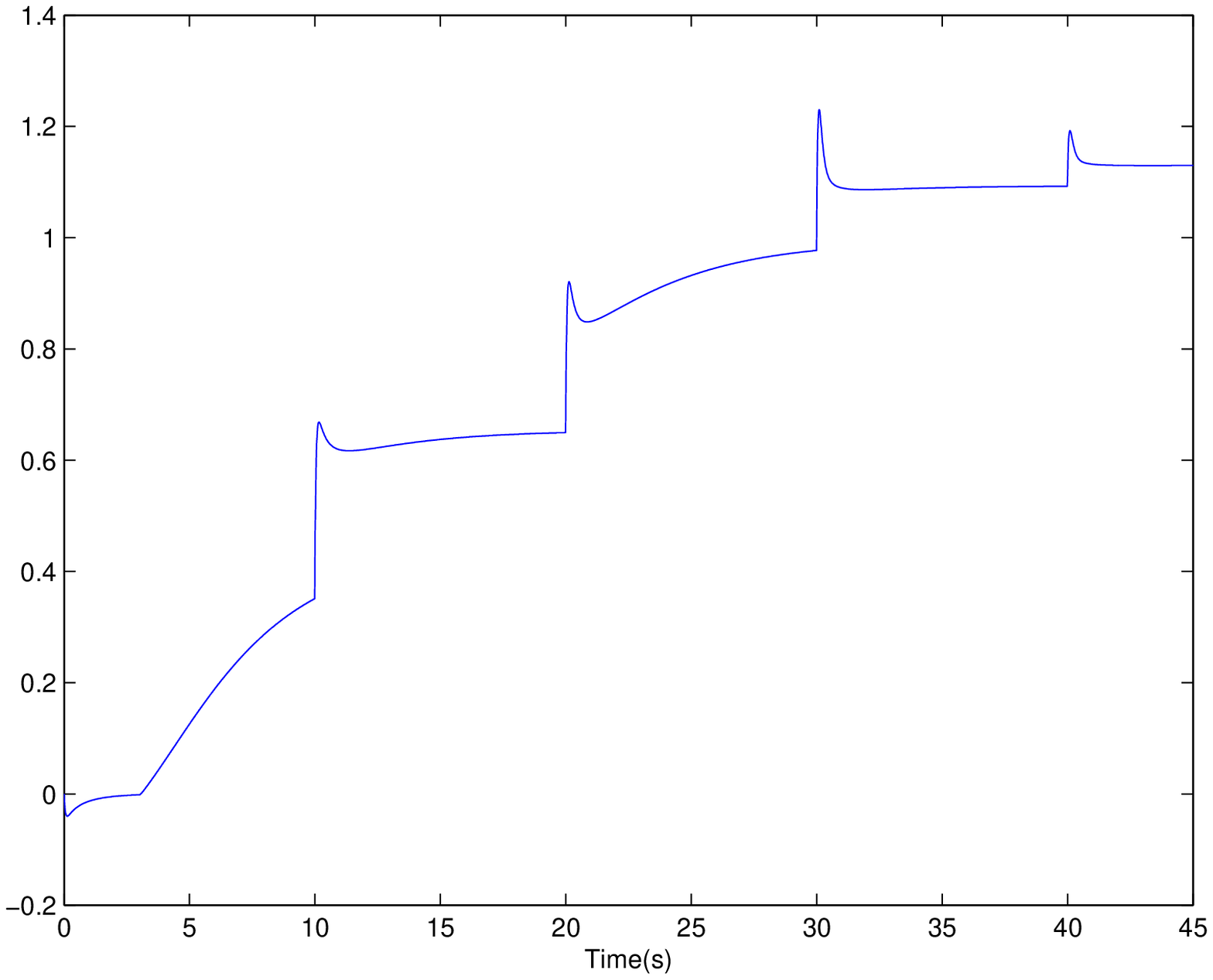}}}%
\subfigure[Setpoint (- .,  black), Reference (- -,  red) and Output (--,  blue)]{
\resizebox*{5.565cm}{!}{\includegraphics{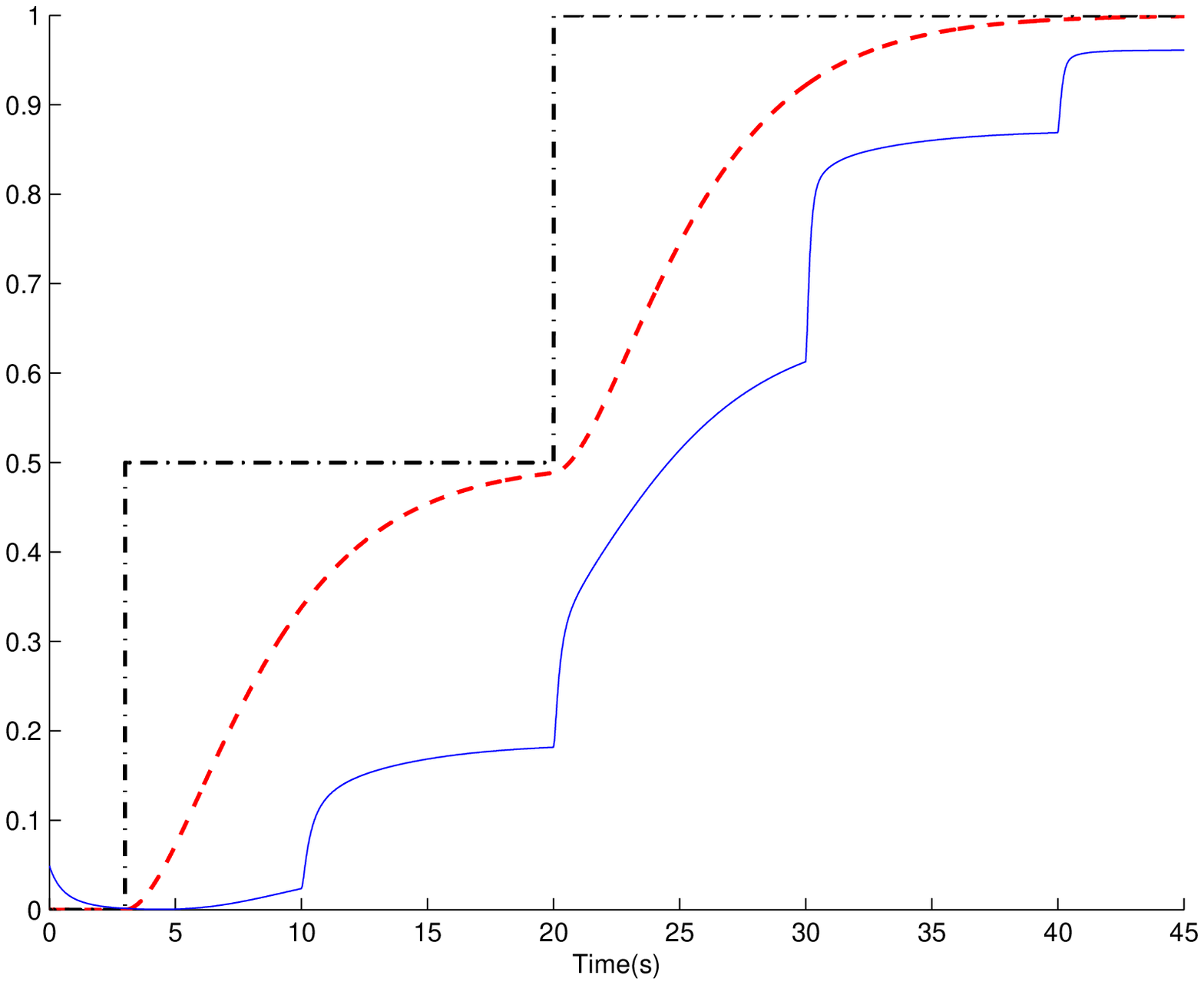}}}
\subfigure[$F$ (- -,  red) and $F_{est}$ (--, blue)]{
\resizebox*{5.565cm}{!}{\includegraphics{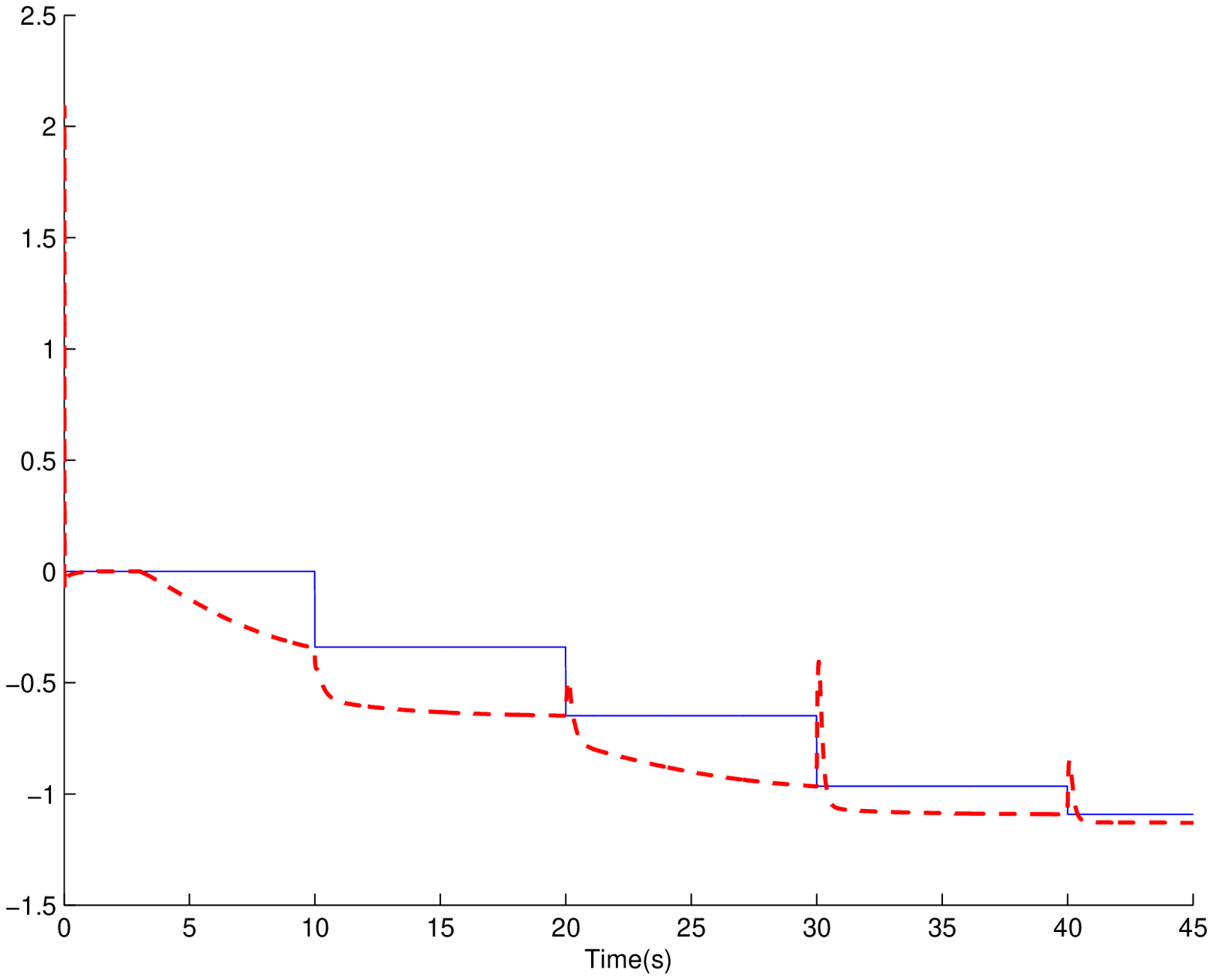}}}%
\caption{$T_{\text{est}}  = 10\text{s}$}
\label{F1000}
\end{center}
\end{figure*}

According to the above control scheme, a ``good'' estimation of $F$ in Equation \eqref{nl1}  plays a key rôle.
Equations \eqref{nl0} and \eqref{nl1} yield the following expression which is used for comparison's sake.
$$F= \frac{3\dot uu^2+2u^3-4 u-\ddot y-3y}{4} $$
Figure \ref{F11} displays excellent results with a sampling time interval $T_{\text{est}}  = 0.01$s for 
estimating $F$.\footnote{$T_{\text{est}}  = 0.01$s is also equal to the sampling time period.} 
The results shown in Figures \ref{F100} and \ref{F1000}
are respectively obtained for $T_{\text{est}}  = 1\text{s}$ and $T_{\text{est}}  = 10\text{s}$. The damages are visible. Figure \ref{F1000} demonstrates that the results with $T_{\text{est}}  = 10\text{s}$ 
cannot be exploited in practice.   
\begin{remark}
Let us emphasize that corrupting noises are neglected here for simplicity's sake. 
\end{remark}

\subsection{A linear academic case}

\subsubsection{Description and control}\label{lac}
\begin{figure*}
\begin{center}
\subfigure[Control]{
\resizebox*{5.565cm}{!}{\includegraphics{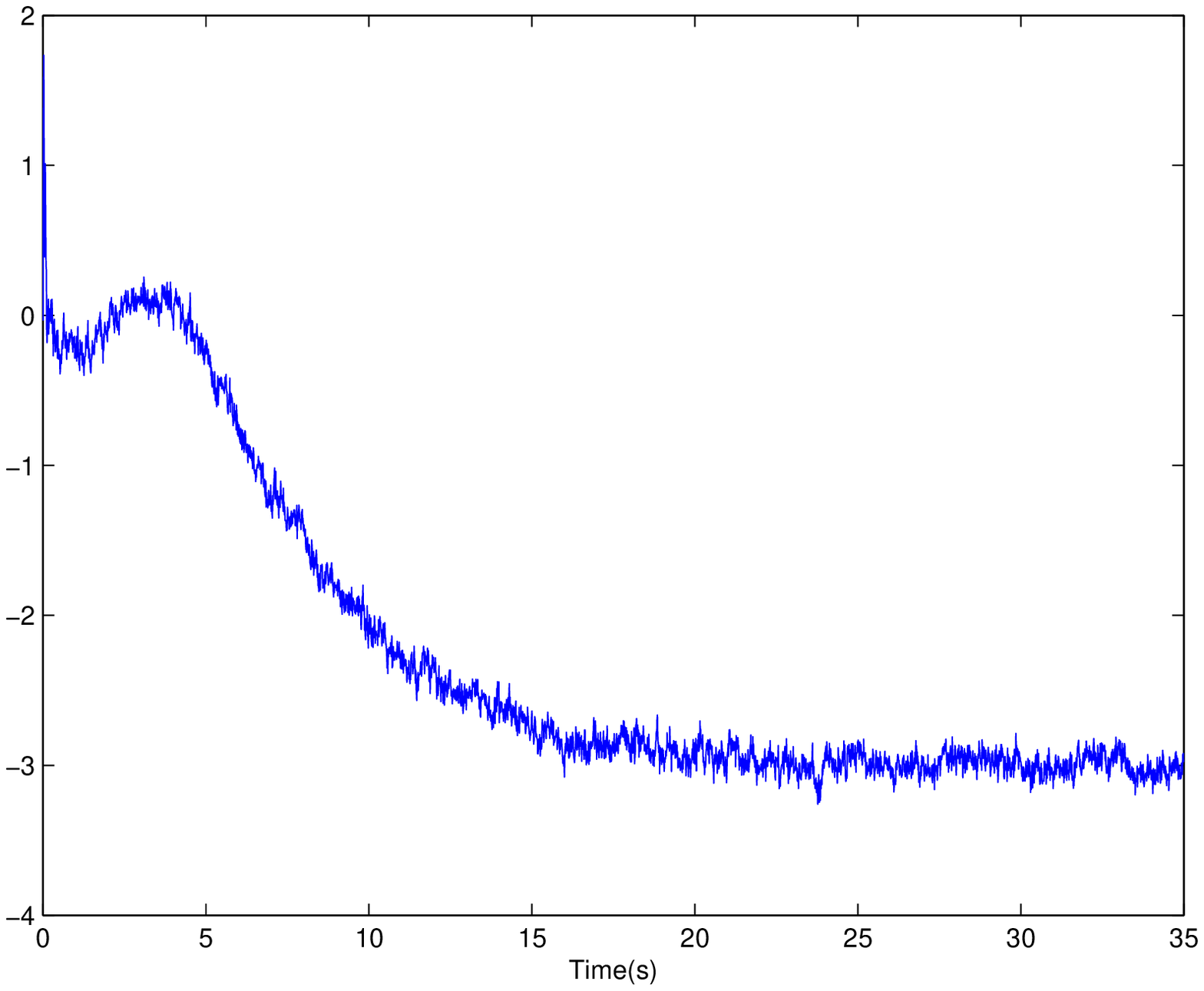}}}%
\subfigure[Setpoint (- .,  black), Reference (- -,  red) and Output (--,  blue)]{
\resizebox*{5.565cm}{!}{\includegraphics{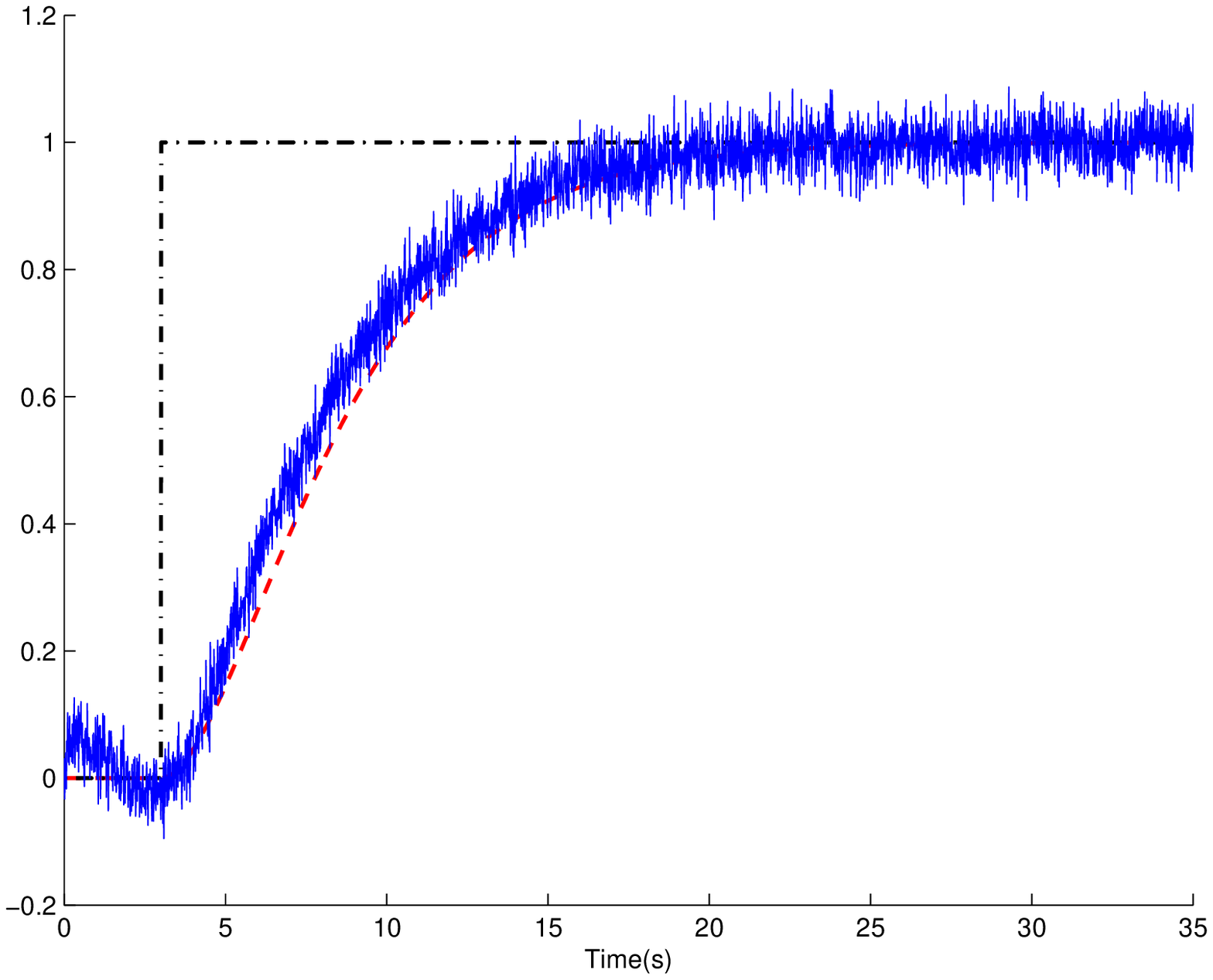}}}%
\caption{Delay free}%
\label{r0}
\end{center}
\end{figure*}

Consider the unstable single-input single-output linear system
\begin{equation}2\dot y-3y=u\label{exemple}\end{equation}
Equation \eqref{nl1} is again used as an ultra-local model. The loop is closed via the iP controller \eqref{fdp}  with some suitable gain $K_P$. As above, in Section \ref{NLDES},
the gain and phase margins are given by Section \ref{IP}.
Stability is therefore ensured with a good robustness.
Figure \ref{r0} displays simulations with $K_P = 1$, a sampling time period $T_e = 0.01\text{s}$, and an additive Gaussian corrupting noise $N(0,0.03)$ on the output. The trajectory tracking is excellent.

\subsubsection{Robustness with respect to a  delayed control}
Introduce a time lag $\tau$ in the control transmission. The transfer function of \eqref{exemple} is no more
$$\frac{1}{2s-3}$$
but
$$\frac{e^{-\tau s}}{2s-3}$$
\begin{remark}
Such delays, which might occur in practice, have already been studied in the literature (see, \textit{e.g.}, \cite{mid,nicu,sipa}).
\end{remark}
\begin{remark}
Systems with transfer functions of the form
$$
T(s)e^{-\tau s}
$$
where $T \in \mathbf{R}(s)$ is a rational function, are according to \cite{retard} the most usual linear delay single-input single-output systems. It is also well known that they are used for approximating 
``complex'' nonlinear systems without delays (see, \textit{e.g.}, \cite{shinskey}). It has been emphasized in \cite{ijc13} that such approximations are becoming useless when applying 
model-free control design. \\
\end{remark}
Assume that we are doing the same computations as in Section \ref{lac}, and, in particular, that $F$ is estimated  with the techniques presented in Section \ref{F}. It amounts saying that we are in fact replacing Equation \eqref{nl1}
by
$$
y(t) = F + u(t - \tau)
$$
The open loop transfer function becomes therefore
$$
T_{1 \tau OLP} = \frac{K_P e^{- \tau s}}{s}
$$
Solving the equation 
$$
T_{1 \tau OLP}(j \omega) = -1
$$
yields 
\begin{equation*}\label{taumax}
\boxed{\tau_{\text{max}}=\frac{\pi}{2K_P}}
\end{equation*} 
\textit{i.e.}, the maximum admissible time lag for stability.

\subsubsection{Computer experiments with delay}
\begin{figure*}
\begin{center}
\subfigure[Control]{
\resizebox*{5.565cm}{!}{\includegraphics{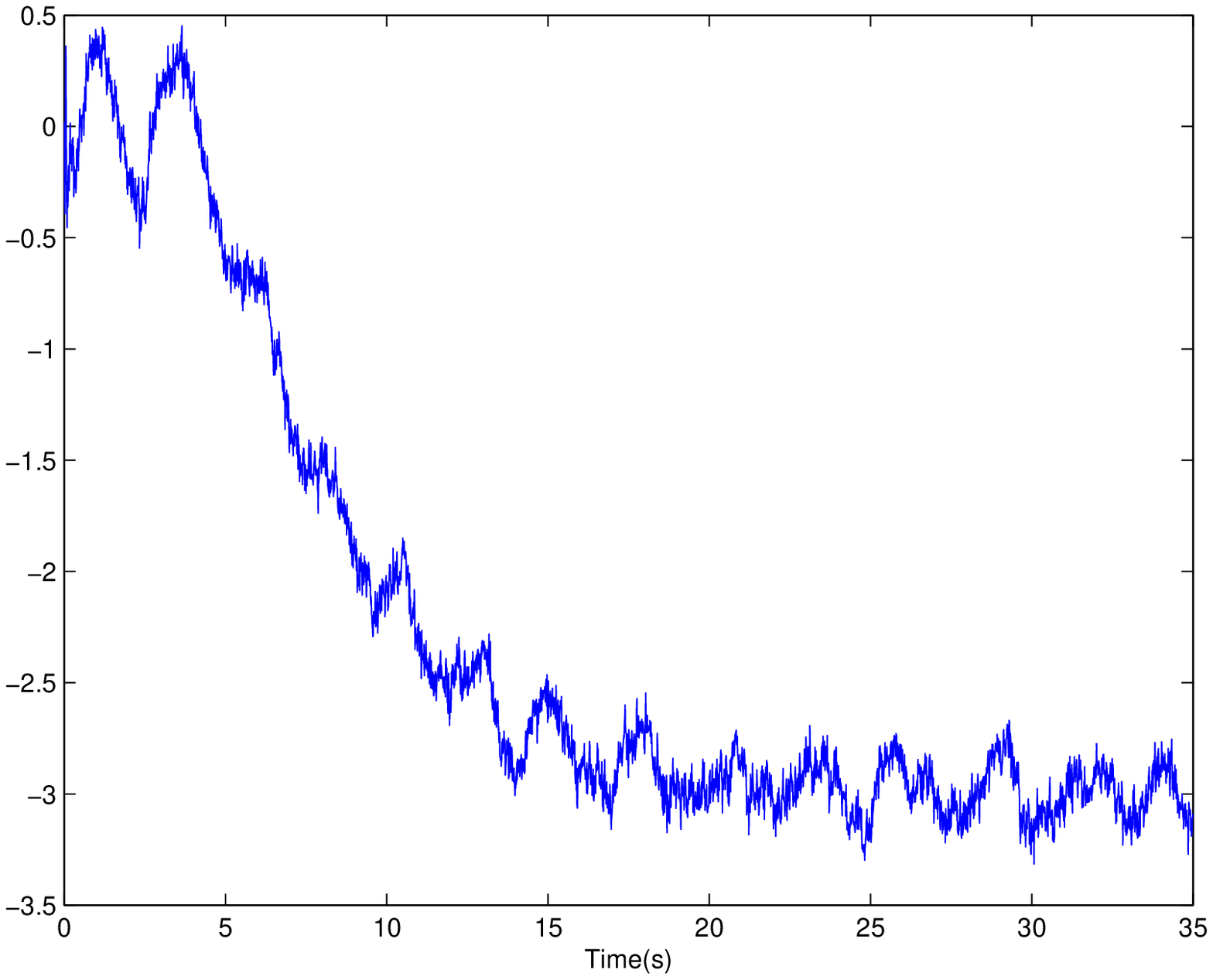}}}%
\subfigure[Setpoint (- .,black), Reference (- -, red) and Output (--, blue)]{
\resizebox*{5.565cm}{!}{\includegraphics{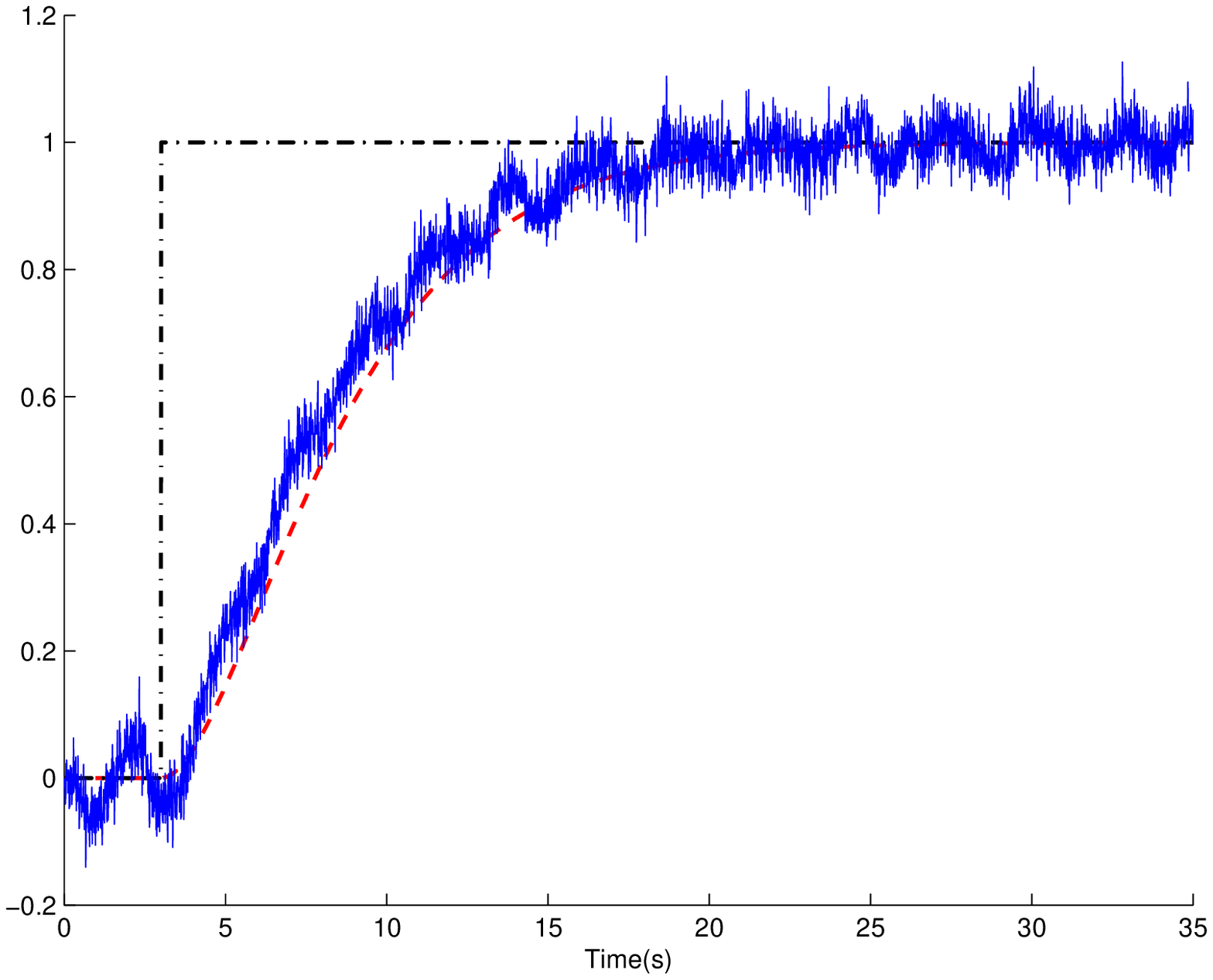}}}%
\caption{Delayed control -- $K_P=1$}%
\label{r1}
\end{center}
\end{figure*}

Figure \ref{r1} displays an excellent stability obtained with a time lag $\tau = 0.2\text{s}$ and $K_P = 1$. Then  $\tau_{\text{max}}\simeq 1.57 \text{s}$.

\begin{figure*}
\begin{center}
\subfigure[Control]{
\resizebox*{5.565cm}{!}{\includegraphics{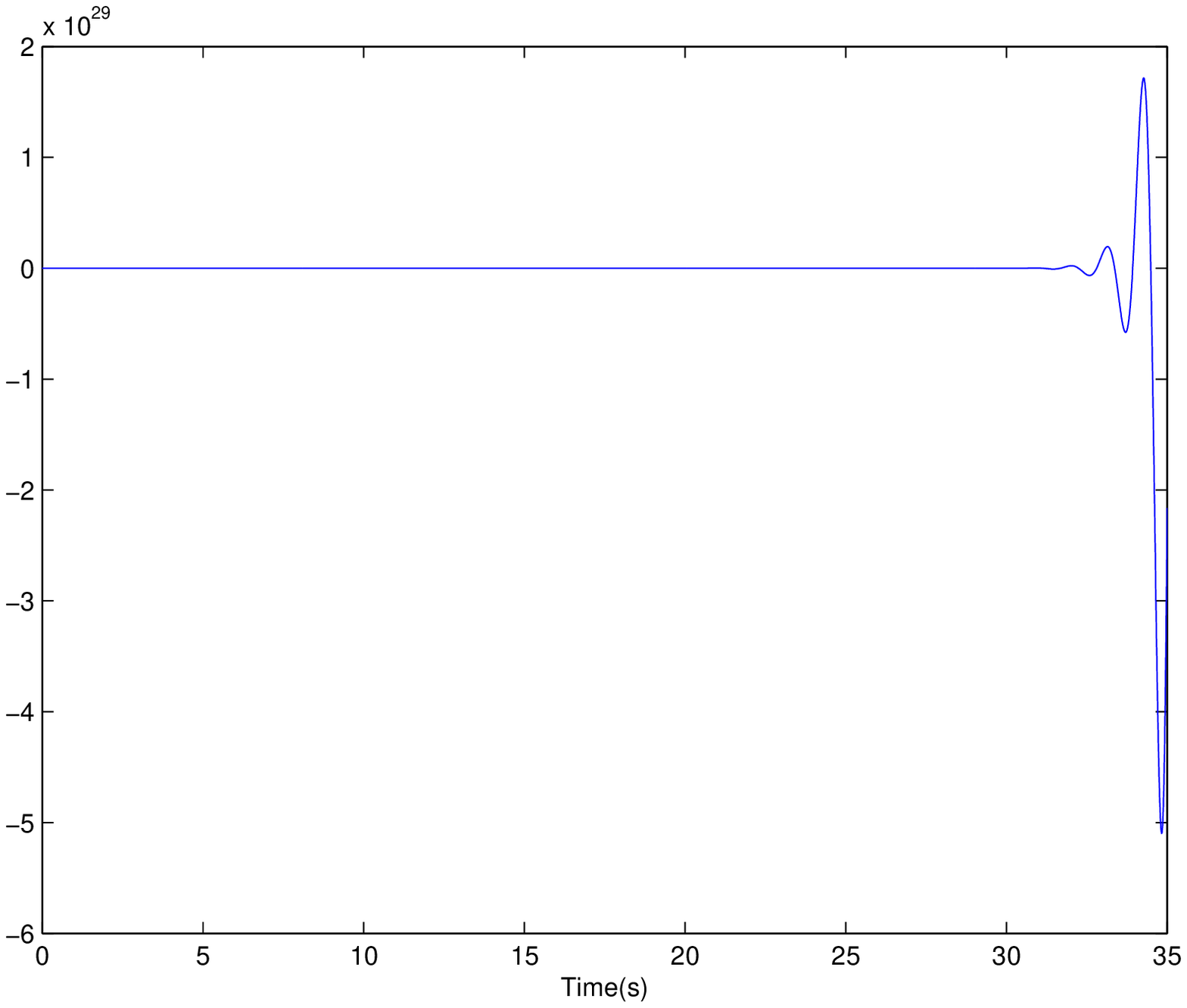}}}%
\subfigure[Setpoint (- .,  black), Reference (- -,  red) and Output (--,  blue)]{
\resizebox*{5.565cm}{!}{\includegraphics{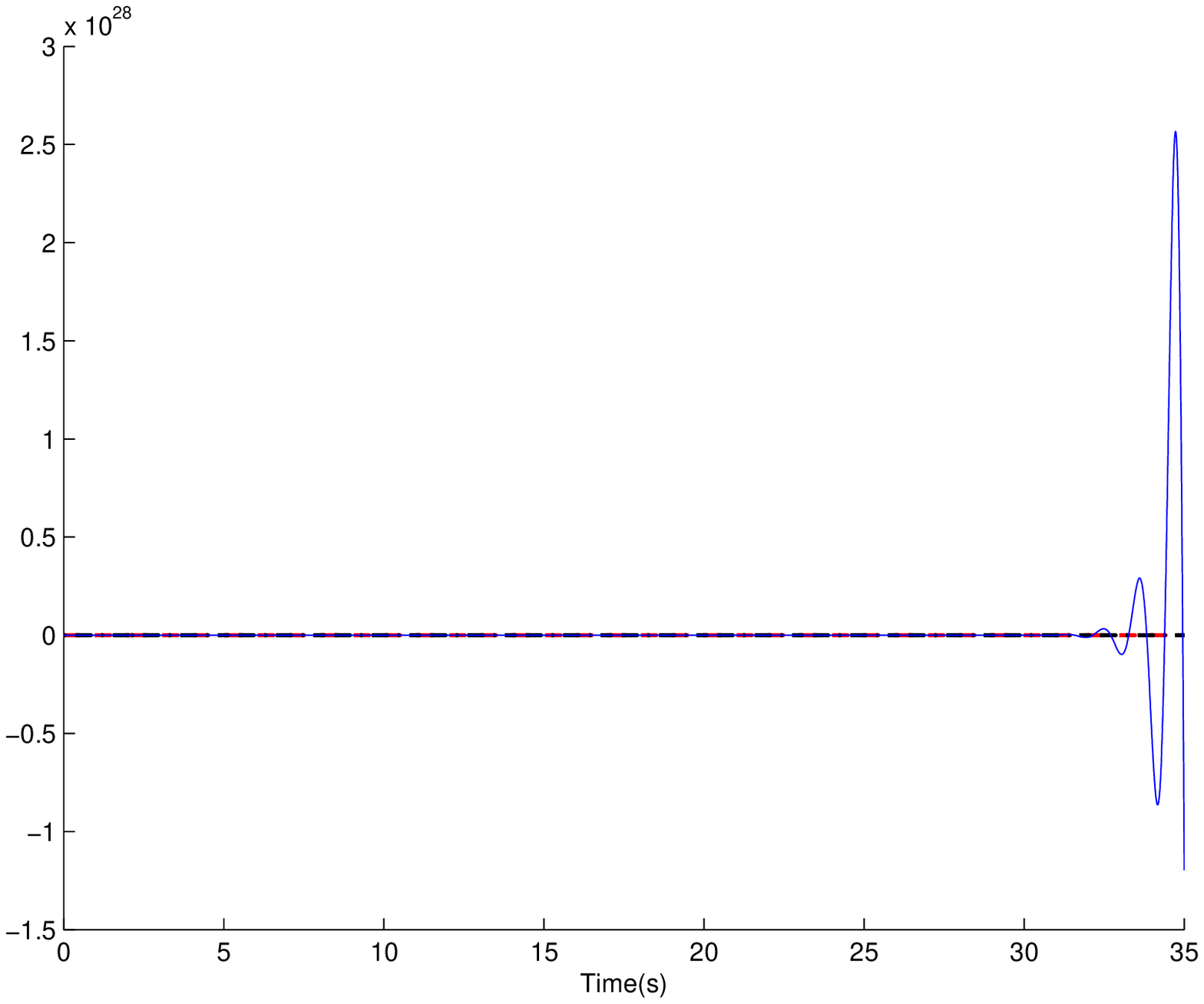}}}%
\caption{Delayed control -- $K_P=10$}%
\label{r2}
\end{center}
\end{figure*}

With the ``high'' gain $K_P=10$, $\tau_{\text{max}}\simeq 0.16 \text{s}$.  Stability is then lost as  shown by  Figure \ref{r2}.


\section{Conclusion}\label{conclusion}
We have demonstrated that the calculations related to stability margins may be easily extended to our recent model-free techniques, where they provide
some new insight on the robustness with respect to delays. As already discussed in \cite{ijc13}, delays, which remain one of the most irritating questions in the model-free setting, do necessitate further 
investigations.\footnote{See also \cite{edf}.}  The key point nevertheless in order to ensure satisfactory performances is in our opinion a ``good'' estimate of $F$. This question, which 
\begin{itemize}
\item has been summarized in Section \ref{F},\footnote{The introduction lists the references of many successful concrete examples.}
\item might become difficult with very severe corrupting noises and/or a poor time sampling,
\end{itemize}
seems unfortunately to be far apart from the stability margins techniques.

\end{document}